\documentclass{amsart}

\usepackage{amssymb}
\usepackage{mathabx}
\usepackage{tikz}
\usepackage{ytableau}
\usepackage{mathptmx}     
\usepackage{helvet}   
\usepackage{courier}
\usepackage{amsfonts} 
\usepackage{rotating}

\usepackage{mathtools}

\newtheorem{theorem}{Theorem}[section]

\theoremstyle{definition}

\theoremstyle{proposition}

\theoremstyle{corollary}

\numberwithin{equation}{section}

\begin{document}

\title[]%
{A visual approach to symmetric chain \\ decompositions of finite Young lattices} 


%    author one information
\author{Terrance Coggins, Robert W. Donley, Jr., Ammara Gondal, and Arnav Krishna}

\address{York College (CUNY), 94-20 Guy R. Brewer Blvd, Jamaica, NY 11451}
\email{terrance.coggins@yorkmail.cuny.edu}

\address{Department of Mathematics and Computer Science,  Queensborough Community College (CUNY), Bayside, NY 11364}
\email{RDonley@qcc.cuny.edu}

\address{Borough of Manhattan Community College (CUNY), 199 Chambers Street, New York, NY 1007}
\email{ammara.gondal001@stu.bmcc.cuny.edu}

\address{Harvard University, Massachusetts Hall, Cambridge, MA 02138}
\email{akrishna@college.harvard.edu}

%    author two information

%\subjclass[2020]{}
%    The 2010 edition of the Mathematics Subject Classification is
%    now available.  If you are citing a classification from the
%    new scheme, use the following input coding instead.
%\keywords{}

\begin{abstract}  
The finite Young lattice $L(m, n)$ is rank-symmetric, rank-unimodal, and has the strong Sperner property.  R. Stanley further conjectured that $L(m, n)$ admits a symmetric chain order.  We show that the order structure on $L(m, n)$ is equivalent to a natural ordering on the lattice points of a dilated $n$-simplex, which in turn corresponds to a weight diagram for the root system of type $A_n$. Lindstr{\" o}m's symmetric chain decompositions for $L(3, n)$ are described completely through pictures.
\end{abstract} 

\maketitle

%section 1
\section{Introduction}

A long-standing conjecture of R. Stanley \cite{StW} concerns the existence of a symmetric chain order on the partially ordered set $L(m, n)$ of  partitions with at most $m$ parts, each of size at most $n$. Several algorithms for the families of $L(3, n)$ and $L(4, n)$ are known, and recently the case of $m=5$ has been settled using computer verification (\cite{bL}, \cite{OH1}, \cite{Rie}, \cite{W1}, \cite{W2}, \cite{Ws}, \cite{Zei}). These approaches are somewhat algorithmic and with few pictures, with \cite{OSZ} as a noteworthy exception.  In this work, we consider the problem as a pen-and-paper exercise on Hasse diagrams.  Our main result is 
\vspace{10pt}

\noindent {\bf Theorem 1.1}:   The order structure on $L(m ,n)$ is equivalent to a natural ordering of the lattice points of an $n$-simplex, dilated by a factor of $m$. In turn, the ordering on this simplex corresponds to a weight diagram of an irreducible finite-dimensional representation for type $A_n.$

\vspace{5pt}

The contents are arranged by section as follows:
\begin{itemize}
\item[(S2)] The finite Young lattices $L(m, n)$,
\item[(S3)] Two views of a simplex,
\item[(S4)] Root systems and weight strings,
\item[(S5)] Symmetric chain decompositions, and
\item[(S6)] Lindstr{\" o}m's algorithm for $L(m, 3).$
\end{itemize}
\vspace{5pt}

During Spring 2024, the first three authors were supported by the Queens Experiences in Discrete Mathematics (QED) REU program, funded by the National Science Foundation, Award Number DMS 2150251. During Summer 2024, the last three authors were supported by the Recruitment and Mentoring in Mathematics (RAMMP) REU program, funded by the National Science Foundation, Award Number DMS 1820731.

% section 2
\section{The Finite Young lattices $L(m, n)$}

\noindent {\bf Definition 2.1:}  Let $k$ be a positive integer.  An {\bf (integer) partition} $\lambda$ of $k$ into $m$ parts is a set of $m$ positive integers $\lambda_i$ such that
$$\lambda_1 + \dots + \lambda_m = k.$$
We write $\lambda \vdash k$ and $|\lambda|=k.$
\vspace{5pt}

It will be convenient to list the parts in non-increasing order, and in examples we simply list the parts as a string of integers.  We also define $\emptyset$ to be the partition of $0$ with no parts.

\vspace{5pt}

\noindent {\bf Example:}  The partitions of 5 with at most 3 parts are
$$5\quad 41 \quad 32 \quad 311 \quad 221.$$

It will be useful to have the corresponding concept with ordered parts.

\vspace{5pt}

\noindent {\bf Definition 2.2:}  Let $k$ be a non-negative integer. A {\bf weak composition} $\alpha$ of $k$ into $n+1$ parts is an ordered set of $n+1$ non-negative  integers $\alpha_i$ such that 
$$\alpha_1 + \dots + \alpha_{n+1} = k.$$

\noindent {\bf Definition 2.3:} The weak compositions of 2 with 3 parts are
$$200\quad 020 \quad 002 \quad 110 \quad 101 \quad 011.$$

A useful visual representation for partitions is given by Young diagrams.

\vspace{5pt}

\noindent {\bf Definition 2.4:}  To construct the {\bf Young diagram} (or {\bf Ferrers diagram}) associated to the partition $\lambda$, we assign a row of squares to each part, listed in non-increasing order and justified to the left.

\vspace{5pt}

\noindent {\bf Example:}  To the partitions listed above, we obtain the following Young diagrams, respectively.

{\tiny $$ 
\begin{ytableau}
\ & \ & \ & \ & \  
\end{ytableau}\qquad\ \ 
\begin{ytableau}
\ & \ & \ & \ \\
\ 
\end{ytableau}\qquad\ \ 
\begin{ytableau}
\ & \ & \ \\
\ & \  
\end{ytableau}\qquad\ \ 
\begin{ytableau}
\ & \  & \ \\
\   \\
\ 
\end{ytableau}\qquad\ \
\begin{ytableau}
\ & \ \\
\ & \  \\
\ 
\end{ytableau}
$$}

\vspace{5pt}

It will be useful to consider sets of partitions with uniform length. When necessary, we extend the partition notation by including zeros as parts.
\vspace{5pt}

\noindent {\bf Definition 2.5:}  Fix positive integers $m, n$. The partially ordered set $L(m, n)$ consists of all partitions with at most $m$ parts, each of size at most $n$. The partial ordering is by entry-wise comparison; that is, $\lambda \le \mu$ if and only if $\lambda_i \le \mu_i$ for all $i$.  Alternatively, $\lambda \le \mu$ if and only if the Young diagram for $\lambda$ fits inside the Young diagram for $\mu$.  

\vspace{5pt}

It will be convenient to define $L(m, n)=\emptyset$ if $m$ or $n$ equals zero.

\vspace{5pt}

\noindent {\bf Definition 2.6:}  We say $\mu$ covers $\lambda$ ($\mu \gtrdot \lambda$) if  $\mu \ge \lambda$ and there is no $\gamma$ with $\mu > \gamma > \lambda$.
\vspace{5pt}

For a covering in $L(m, n)$, there is exactly one $i$ such that $\mu_i = \lambda_i + 1$. Alternatively, the Young diagram for $\lambda$ is obtained by removing one square from the Young diagram for $\mu$.

\vspace{5pt}

To construct the Hasse diagram of $L(m, n)$, the maximal element is a string of $m$ values of $n$; the corresponding Young diagram is an $m\times n$ rectangle.  Links are determined by covering, and levels correspond to partitions with the same sum of parts or Young diagrams with the same number of squares.

\vspace{5pt}

\noindent {\bf Example:}  Consider the ordering for  $L(4, 3)$.  To see that $22 \le 32$, we compare $2200$ and $3200$ entry-wise.  On the other hand, $1111$ and $22$ are not comparable.

\vspace{5pt}

\noindent {\bf Example:}  When $m=1$, $L(m, n)$ is a chain of length $n$, and a similar result holds for $n=1.$ 

\vspace{5pt}

\noindent {\bf Example:}  When $n=2$, we obtain triangular shaped diagrams as seen in figure 1.

\begin{figure}
\centering
\includegraphics[width=.12\textwidth]{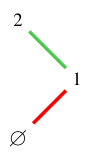}\qquad\qquad
\includegraphics[width=.2\textwidth]{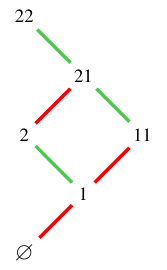}\qquad
\includegraphics[width=.32\textwidth]{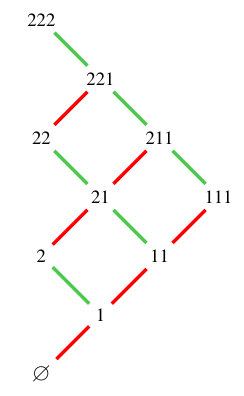}
\caption{Finite Young lattices $L(1, 2)$, $L(2, 2)$, $L(3, 2)$}
\label{}
\end{figure}

\vspace{5pt}

\noindent {\bf Defintion 2.7:} The {\bf conjugate} of the partition $\lambda,$ denoted by $\lambda',$  is the partition with parts given by number of squares in each column of the Young diagram for $\lambda.$ The Young diagram for $\lambda'$ is similar to a matrix transpose, which interchanges rows and columns.

\vspace{5pt}

For the partitions of 5 above, the conjugates are given by 11111, 21111, 221, 311, and 32. The partition 311 is self-conjugate.

As partially ordered sets, $L(m, n)$ and $L(n, m)$ are equivalent through conjugation and have similar Hasse diagrams.  Compare the Hasse diagrams for $L(3, 2)$ and $L(2, 3)$ in figures 1 and 3.

\vspace{5pt} 

Finally, there are $\begin{pmatrix} m+n \\ m \end{pmatrix} = \frac{(m+n)!}{m!n!}$ elements in $L(m, n)$.  The rank numbers for $L(m, n)$ are given by the $q$-binomial coefficient
$$\begin{bmatrix} m+n \\ m\end{bmatrix}_q = \frac{\bf(m+n)!}{\bf (m)!\ (n)!},$$
where 
$${\bf (m)!} = (1+q) (1+ q + q^2)\dots (1+ q + \dots + q^{m-1}) =\frac{(1-q)(1-q^2)\dots (1-q^m)}{(1-q)^m}.$$
The coefficient of $q^k$ in $\begin{bmatrix} m+n \\ m\end{bmatrix}_q$ gives the number of partitions of $k$ in $L(m, n)$.

\noindent {\bf Example:}  As seen in figure 5 below, $L(3, 3)$ has $\begin{pmatrix} 6 \\ 3 \end{pmatrix} = 20$ elements.  The corresponding rank sequence is given by the coefficients of
\begin{eqnarray}\begin{bmatrix} 6 \\ 3\end{bmatrix}_q &= &\frac{(1-q^6)(1-q^5)(1-q^4)(1-q^3)(1-q^2)(1-q)}{(1-q^3)^2(1-q^2)^2(1-q)^2}\notag\\
\ &= &1+q + 2q^2 + 3q^3 + 3q^4 + 3q^5 + 3q^6 + 2q^7 + q^8 + q^9,\notag
\end{eqnarray}
also verified by this figure.

Several interesting identities for $q$-binomial coefficients are expressed through the Hasse diagram for $L(m, n).$ For instance, consider the the identities
$$\begin{bmatrix} m+n \\ m\end{bmatrix}_q = q^n\begin{bmatrix} m+n-1 \\ m-1\end{bmatrix}_q + \begin{bmatrix} m+n-1 \\ m\end{bmatrix}_q  = \begin{bmatrix} m+n-1 \\ m-1\end{bmatrix}_q + q^m\begin{bmatrix} m+n-1 \\ m\end{bmatrix}_q.$$
The first equality represents the splitting of the Hasse diagram for $L(m, n)$  into subdiagrams for $L(m-1, n)$ and $L(m, n-1)$, with the second diagram shifted up by $n$ levels. This equality simply counts the partitions with and without a part of size $n$. The second equality has a similar interpretation in terms of number of parts.

\vspace{5pt}

{\bf Example:} In figure 1, to construct $L(3, 2)$ from $L(2, 2)$, we append a 2 to each partition in $L(2, 2)$ and attach these partitions one level above the chain $L(3, 1)$.  For the second equality, we obtain $L(3, 2)$ from $L(2, 2)$ in a similar manner by adding a chain to the top of the diagram.
\vspace{5pt}

{\bf Example:} Figure 6 shows several decompostions for the diagram of $L(3, 3).$ The first diagram splits into subdiagrams with or with a part of size 3.  Furthemore, we iterate the first equality for $L(3, 3)$ to obtain the decomposition 
$$L(3,3) = L(2, 3) \cup L(2, 2) \cup L(2, 1) \cup L(2, 0)$$
according to largest parts.  The third decomposition splits the diagram by number of parts.
\vspace{5pt}

See \cite{And} and \cite{AnE} for the general theory of integer partitions.  For order theory and integer partitions, see also \cite{BSa}, \cite{StA}, and \cite{StE}.

%section 3
\section{Two views of a simplex}

While the Hasse diagram for $L(m, n)$ becomes unwieldy when arranged by levels, simplices offer a natural visualization that displays the ordering.

\vspace{5pt}

\noindent {\bf Definition 3.1:}  The $(m-1)$-simplex $\Delta'(m-1)$ is the set of all points $(x_1, \dots, x_m)$ in $\mathbb{R}^{m}$ such that
$$x_i\ge 0, \qquad x_1 + \dots + x_m =1.$$
In other words, $\Delta'(m-1)$ is the convex hull of the standard basis vectors $\{e_1, \dots, e_m\}$ in $\mathbb{R}^m$.

On the other hand, the simplex $\Delta'(m+1) $ is combinatorially equivalent to the polytope $\Delta(m)$ defined by all $(y_1, \dots, y_{m})$
such that $1\ge y_1 \ge \dots \ge y_{m} \ge 0.$

\vspace{5pt}

\noindent {\bf Example:} Compare the two-dimensional simplices in figure 2.

\vspace{5pt}

\noindent {\bf Proposition 3.2:}  $L(m, n)$ is the set of integer points in $n\Delta(m)$,  the simplex in $\mathbb{R}^m$ dilated by a factor of $n.$

\vspace{5pt}

The simplex $\Delta'(m)$ allows another interpretation of $L(m, n)$ through multiplicities of parts.

\vspace{5pt}

\noindent {\bf Definition 3.3:}  The poset $L'(m, n)$ consists of all weak compositions of $m$ with $n+1$ parts.  We say that $\alpha$ covers $\beta$ if $\alpha$ and $\beta$ have the same parts, except that a value of 1 is subtracted from exactly one $\alpha_i$ and added to $\alpha_{i+1}$ to obtain $\beta$.
\vspace{10pt}

\noindent {\bf Proposition 3.4:}  $L'(m, n)$ is the set of integer points in $m\Delta'(n)$,  the simplex in $\mathbb{R}^{n+1}$ dilated by a factor of $m.$

\vspace{5pt}

\begin{figure}
\centering
\includegraphics[width=.3\textwidth]{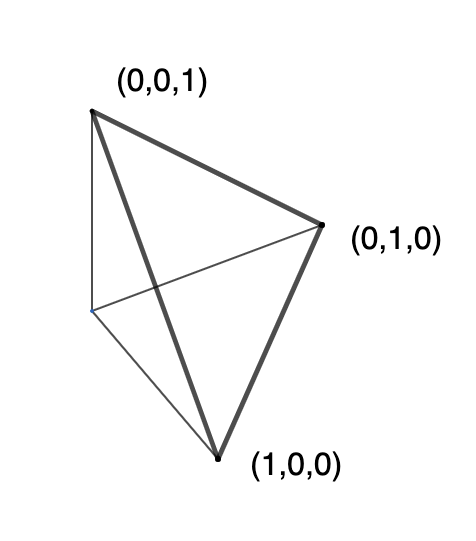}\qquad\qquad
\includegraphics[width=.3\textwidth]{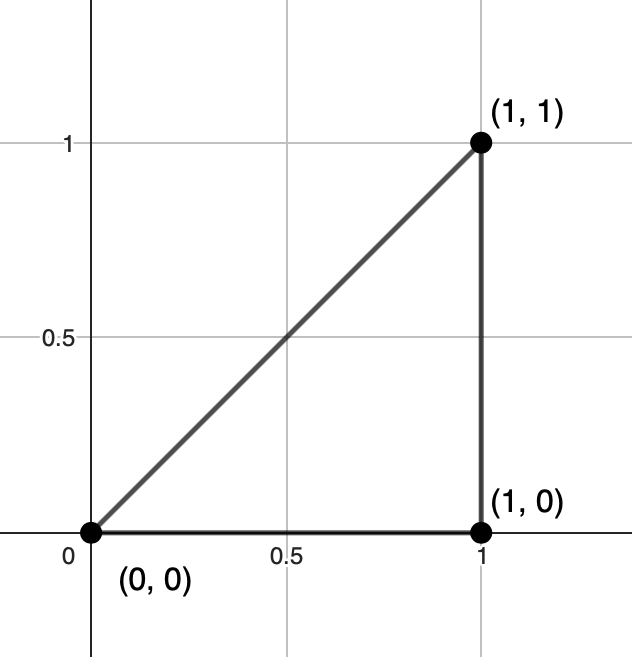}\qquad
\caption{The simplices $\Delta'(2)$ and $\Delta(2)$}
\label{}
\end{figure}

\noindent {\bf Example:}  In $L'(4, 3),$ the composition 1120 covers 0220, 1030, and 1111. Compare figures 10 and  11 for the Hasse diagrams of $L(4, 3)$. and $L'(4, 3).$
\vspace{5pt}

\noindent {\bf Proposition 3.5:}  There is an order-preserving one-one correspondence from $L(m, n)$ to $L'(m, n).$  To implement, we send each partition, with additional zero parts included, to a weak composition that records the multiplicities of each part, listed in decreasing order of part size.

\vspace{5pt}

\noindent {\it Proof.}  The one-one correspondence property is immediately checked.  We check that the covering relations are preserved.  If $\mu$ covers  $\lambda$ in $L(m, n),$ then there is a unique part $x$ that has been decreased by 1 in $\mu$ to obtain $\lambda.$ In other words, the multiplicity of $x$ has decreased by 1 while the multiplicity of $x-1$ has increased by 1. QED

\vspace{5pt}

\noindent {\bf Example:}  In $L(4, 3)$, the part sizes are 3, 2, 1, and 0.  The partition 3211, which corresponds to the weak composition 1120, covers 2211, 3111, and 321. In turn, these partitions map to the compositions 0220, 1030, and 1111.  Here we use 3210 instead of 321 to obtain 1111.

% section 4
\section{Root systems and weight strings}

The ordering on $L'(m, n)$ is a familiar example from root systems and reflection groups in representation theory. See \cite{FH}, \cite{Hu}, or \cite{Kn} for basic theory and notation for root systems and finite-dimensional representations.

\vspace{5pt}

\noindent {\bf Definition 4.1:}  Let $\{e_i\}$ be the standard basis for $\mathbb{R}^{n+1}$. The root system for type $A_{n}$ is the set
$$\Delta = \{e_i-e_j\ |\ 1\le i, j \le n+1, i\ne j\}.$$  We choose the set of positive roots
$$\Delta^+ = \{e_i-e_j\ | 1\le i <j \le n+1 \}$$
with corresponding set of $n$ simple roots
$$\Pi = \{\alpha_i=e_i-e_{i+1}\}.$$
 
Every $\beta$ in $\Delta^+$ can be represented uniquely as a linear combination of simple roots with nonnegative integer coefficients; that is,
$$\beta = \sum\limits_{i=1}^{n} m_i\alpha_i$$
with $m_i\ge 0.$ In turn, $\Pi$ forms a basis $B'$ of the subspace $V$ of $\mathbb{R}^{n+1}$ of dimension $n$ given by
$$V=\{\sum\limits_{i=1}^{n+1} m_i e_i\ \ | \ \sum\limits_{i=1}^{n+1} m_i=0 \}.$$
\vspace{5pt}

\noindent {\bf Definition 4.2:} The {\bf weight lattice} $\Gamma$ is the set of all elements $\omega$ in $V$ such that $\frac{2\langle \omega, \alpha\rangle}{\langle \alpha, \alpha \rangle}$ is an integer for all $\alpha$ in $\Pi.$ The {\bf fundamental weight} $\omega_i$ in $\Gamma$ satisfies $\frac{2\langle \omega_i, \alpha_j\rangle}{\langle \alpha_j, \alpha_j \rangle}=\delta_{ij}$, where $\delta_{ij}=1$ if $i=j$ and 0 otherwise.
\vspace{5pt}

We form a basis $B$ of $\mathbb{R}^{n+1}$ by including $\omega_0=e_1 + \dots + e_{n+1}$ with $B'.$. When working with weight diagrams, it will be convenient to translate $V$ by an appropriate multiple of $\alpha_0$ so that all elements of interest have nonnegative integer entires.   In this case, we represent the fundamental weights by
$$\omega_i = e_1 + \dots + e_i$$ 
for $1\le i < n+1.$

\vspace{5pt}

\noindent {\bf Definition 4.3:}  A {\bf dominant weight} $\lambda$ is an element of the weight lattice such that $\langle \lambda, \alpha_i\rangle \ge 0$ for all simple roots $\alpha_i.$ Alternatively, $\lambda$ is a nonnegative integer linear combination of fundamental weights.   
\vspace{5pt}

As elements in $\mathbb{R}^{n+1}$, dominant weights are of the form $\lambda = \sum\limits_{i=1}^{n+1} m_ie_i$ with integer $m_i$ such that $0\le m_{n+1} \le m_{n} \le \dots \le m_1.$
\vspace{5pt}

\noindent {\bf Definition 4.4:}  Choose a dominant weight $\lambda$, and form the set $S$ of all points in $\mathbb{R}^{n+1}$ whose coordinates are the same as $\lambda$ up to a permutation.  The {\bf weight diagram} $W(\lambda)$ associated to the highest weight $\lambda$ is the set of all weights contained in the convex hull of $S$ in $\mathbb{R}^{n+1}$.  
\vspace{5pt}

Every element $\omega$ of $W(\lambda)$ can be represented in the form $\omega = \lambda - \sum\limits_{i=1}^{n+1} m_i\alpha_i$.  This representation need not be unique.
\vspace{5pt}

\noindent {\bf Example:}  We consider only one type of highest weight: $\lambda = m\omega_1 = me_1.$  In this case, $W(\lambda)$ is the convex hull of the weights
$$\{me_1, me_2, \dots, me_{n+1}\}.$$
That is, $W(\lambda)$ is the set of weights contained in the dilated $n$-simplex
$$m\Delta'(n) = \{(x_1, \dots, x_{n+1})\ |\ x_i\ge 0,\ \ \sum x_i = m\}.$$

\noindent {\bf Example:}  Let $U=\mathcal{P}_{m, n+1}$ be the set of homogenous polynomials of degree $m$ in $n+1$ complex variables $\{z_i\}$.  A basis for $U$ is given by the monomials
$$\{z_1^{m_1}\dots z_{n+1}^{m_{n+1}}\}.$$  These indices $(m_1, \dots, m_{n+1})$ form the set of weak compositions of $m$ with $n+1$ parts, so
$$\dim_\mathbb{C} (U) = \begin{pmatrix} m+n \\ m\end{pmatrix}.$$
The group $G=GL(n+1, \mathbb{C})$ acts on $U$ by
$$[\pi_{m}(g)P](z) = P(g^Tz),$$
where $g^T$ denotes the matrix transpose of $g$ in $G.$

The diagonal matrix $diag(t_1, \dots, t_{n+1})$ acts on the monomial $P(z) = z_1^{m_1}\dots z_{n+1}^{m_{n+1}}$ by the scalar $t_1^{m_1}\dots t_{n+1}^{m_{n+1}}$. We associate to $P(z)$ the weight $(m_1, \dots, m_{n+1})$, for which $m_1+\dots+ m_{n+1} = m.$  This representation has highest weight $\lambda=me_1$, each weight has multiplicity 1, and  $W(\lambda) = m\Delta'(n).$

\vspace{5pt}

With these identifications, we have a both a natural visualization and consistent  language to describe the order structure on $L'(m, n).$ 

\vspace{5pt}

\noindent {\bf Definition 4.5:}  The weight string through $\gamma$ along the root $\alpha$  in $W(\lambda)$ is a maximal subset of weights in $W(\lambda)$ of the form
$$\dots,\ \ \gamma - 2\alpha,\ \ \gamma-\alpha,\ \ \gamma, \ \ \gamma+\alpha,\ \  \gamma+2\alpha, \ \ \gamma + 3\alpha, \dots.$$

\vspace{10pt}

\noindent {\bf Definition 4.6:}. To color the edges of the Hasse diagram for  $L'(m, n),$ we first choose colors for each simple root $\alpha_i$.  Then the edges for each weight string along $\alpha_i$ are given the same color.  We color $L(m, n)$ by equivalence.

\vspace{10pt}

\noindent {\bf Example (Colors for $L(n, 2)$):}  We use green for $\alpha_1=e_1-e_2$ and red for $\alpha_2=e_2-e_3.$   In the partition notation of figure 1, we use green when we reduce a 2 to 1 and red when we remove a 1. 

\vspace{10pt}

\noindent {\bf Example (Colors for $L(n, 3)$):}  In figures 3 and 5 through 13, we use green for $\alpha_1=e_1-e_2$, red for $\alpha_2=e_2-e_3,$ and blue for $\alpha_3=e_3-e_4.$   For Hasse diagrams labeled with Young diagrams or partition notation, we use green when we reduce a 3 to 2, red when we remove reduce a 2 to 1, and blue when we remove a 1.  The rules in $L'(n, 3)$ apply when we shift a 1 between the corresponding adjacent entries.

\vspace{10pt}

\noindent {\bf Example (Figures 10 and 11):}. Consider the red weight string through $\lambda=3222$ along $\alpha_2=e_2-e_3$ with Young diagrams, partitions, and weak compositions, respectively.  

{$$\tiny
\begin{ytableau}
\ & \ & \ \\
\ & \  \\
\ & \  \\
\ & \ 
\end{ytableau}\qquad\ {\small \gtrdot} \qquad
\begin{ytableau}
\ & \ & \ \\
\ & \  \\
\ & \ \\
\
\end{ytableau}\qquad \gtrdot \qquad
\begin{ytableau}
\ & \ & \  \\
\ & \ \\
\ \\
\
\end{ytableau}\qquad\ \gtrdot \qquad 
\begin{ytableau}
\ & \ & \  \\
\ \\ 
\ \\
\
\end{ytableau}
$$}

$$3222\ \ \gtrdot \ \ 3221 \ \ \gtrdot \ \ 3211 \ \ \gtrdot \ \ 3111$$
$$1300\ \ \gtrdot \ \ 1210 \ \ \gtrdot \ \ 1120 \ \ \gtrdot \ \ 1030$$
\vspace{5pt}

Finally, we note that, while conjugation gives an equivalence of orderings between $L(m, n)$ to $L(n, m)$, the colored Hasse diagrams need not agree.  For instance, compare the diagram for $L(3, 2)$ in figure 1 with $L(2, 3)$ in figure 3.  As noted above, the number of colors equals the number of simple roots.

\begin{figure}
\centering
\includegraphics[width=.32\textwidth]{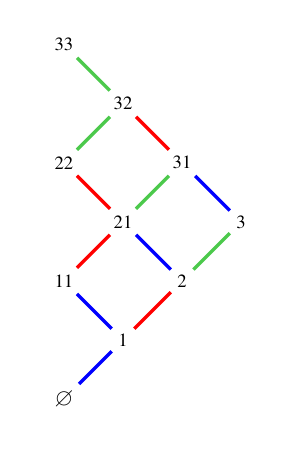}\qquad
\includegraphics[width=.4\textwidth]{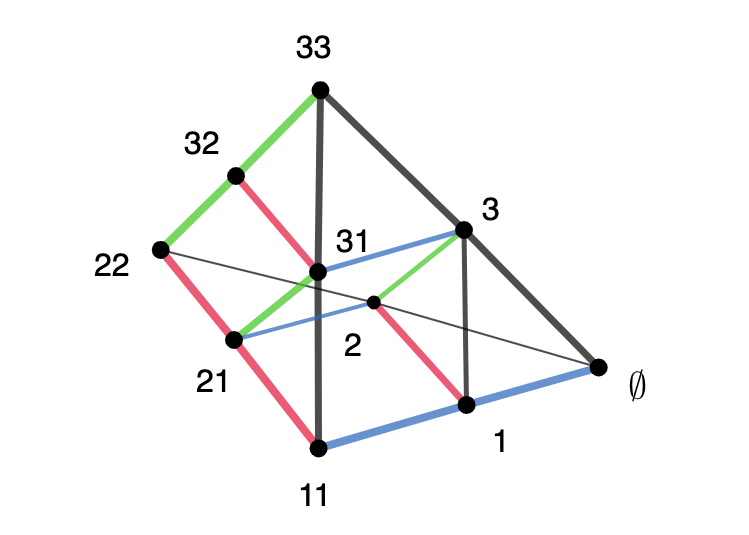}
\caption{Finite Young lattice $L(2, 3)$ with three simple roots}
\label{}
\end{figure}

%section 5
\section{Symmetric chain decompositions}

We recall some definitions needed to clarify the rank function on $L(m, n)$ and define a symmetric chain order for a partially ordered set $(P, \le)$.
\vspace{5pt}

\noindent {\bf Definition 5.1:}. A {\bf saturated chain} in a partially ordered set $P$ is a sequence of consecutive coverings.  In the Hasse diagram for $P$, a saturated chain is a vertical path with no gaps.The length of the chain is the number of coverings in the sequence.
\vspace{5pt}

Examples of saturated chains are given by the weight strings in section 4.
\vspace{5pt}

\noindent {\bf Definition 5.2:}. A finite partially ordered set $P$ is called {\bf graded} with ${\hat 0}$ and ${\hat 1}$ if it has unique minimum and maximum elements ${\hat 0}$ and ${\hat 1}$, and, for each $x$ in $P$, every saturated chain from the minimum element ${\hat 0}$ to $x$ has the same length.  For each $x$ in $P$, the rank function $\rho(x)$ is the length of any saturated chain to $x$.  The {\bf level} for rank $s$ is the set of all elements in $P$ with rank equal to $s.$. The height $ht(P)$ of $P$ is the largest rank.
\vspace{5pt}

For $L(m, n)$, $\rho(\lambda)$ is either the number of squares in the Young diagram for $\lambda$ or the sum of the parts of $\lambda$ as a partition.  If $\lambda$ is represented as a weak composition $\alpha$, $\rho(x)=\sum i\alpha_i.$
\vspace{5pt}

Next, $L(m, n)$ is {\bf self-dual}; that is, there exists an involution of $L(m, n)$ that reverses the partial order.  The natural involution is given by  complementation.

\vspace{5pt}

\noindent {\bf Definition 5.3:}. The {\bf complement} of $\lambda$ in $L(m, n)$ equals the partition $\mu=\lambda^*$ with parts $n-\lambda_i$.  To list as a string of non-increasing parts, we reverse the order also.  The rank of $\lambda^*$ in $L(m, n)$  is $mn-|\lambda |.$
\vspace{5pt}

To obtain the complement of $\lambda$ in $L(m, n)$ as a Young diagram, we remove the squares that correspond to $\lambda$ in the $m\times n$ rectangle and rotate the remaining diagram by 180 degrees.

\vspace{5pt}
{\bf Example:}. Consider the partition 322 in $L(4, 3).$  The complement has parts 1, 2, 2, with corresponding Young diagram
$$\begin{ytableau}
\ & \ & \  \\
\ & \  \\
\ & \ 
\end{ytableau} \quad \to \quad 
\begin{ytableau}
 X & X & X  & \ \\
X & X & \ & \ \\ 
X & X & \ & \ \\
\end{ytableau}\quad \to\quad
\begin{ytableau}
\ & \  \\
\ & \ \\ 
\ 
\end{ytableau}.
$$ 

\noindent {\bf Definition 5.4:}. A {\bf symmetric chain} is a saturated chain $x_1 \le x_2 \le \dots \le x_k$  such that $\rho(x_1) + \rho(x_k) = ht(P).$ That is, the chain begins and ends on opposite levels. 
\vspace{5pt}

For example, a symmetric chain in $L(m, n)$ that begins in rank $s$ ends in rank $mn-s.$
\vspace{5pt}

\noindent {\bf Definition 5.5:}. A {\bf symmetric chain decomposition} of $P$ is a collection of symmetric chains such that each element of $P$ is contained in exactly one such chain.  A partially ordered set that admits a symmetric chain decomposition is called a {\bf symmetric chain order.}
\vspace{5pt}

\noindent {\bf Conjecture (Stanley \cite{StW}):}. $L(m, n)$ admits a symmetric chain order. 
\vspace{5pt}

As noted in the introduction, it is known that $L(m, n)$ is a symmetric chain order for $m\le 5$.

\vspace{5pt}

\noindent {\bf Example ($L(m, 2)$):}. In figure 1, we simply alternate between green and red edges, starting from the left.  When $m$ is even, there is a singleton chain. When $m$ is odd, the smallest chain has length 2.
See the leftmost face of the simplices in figures 3, 8, and 9. 

%section 6
\section{ Lindstr{\" o}m's algorithm for $L(m, 3)$}

We use the preceding discussion to illustrate Lindstr{\"o}m's algorithm for $L(3, n)$ \cite{bL}.  In fact, we given the ordering on the equivalent $L(m, 3)$; all weak compositions of interest then have length 4. The algorithm splits into odd and even cases, recursively defined with periods 2 and 4, respectively.  In each case, we fit a simplex with a symmetric chain decomposition symmetrically into a larger simplex and fill the remainder with symmetric chains, as indicated in figure 4.

\begin{figure}
\centering
\includegraphics[width=.4\textwidth]{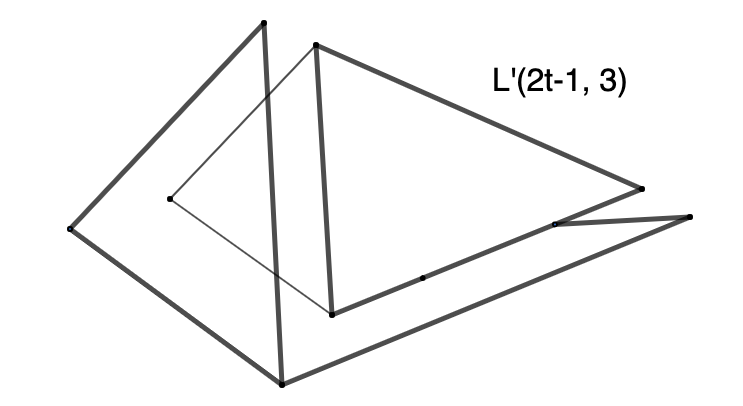}\qquad
\includegraphics[width=.43\textwidth]{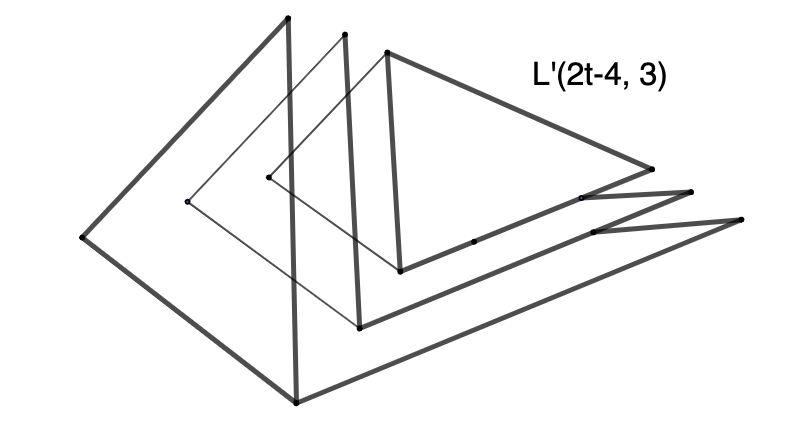}
\caption{Finite Young lattices $L'(2t+1, 3)$ and $L'(2t, 3)$}
\label{}
\end{figure}

Symmetric chain decompositions for $L'(2t+1, 3)$  are directly defined.  We work with weak compositions instead of partitions. First, $L'(1, 3)$ is a chain with symmetric chain decomposition $1000 \gtrdot 0100 \gtrdot 0010 \gtrdot 0001.$  Suppose the symmetric chain decomposition exists for $L'(2t-1, 3).$  We include these chains symmetrically in $L'(2t+1, 3)$ by sending each weak composition $abcd$ to $(a+1)bc(d+1).$  For each partition, we append a part of size 3.   The remaining chains in $L'(2t+1, 3)$ consist of weak compositions with first or last entry equal to 0.  

Each remaining symmetric chain starts on the edge $ab00$ and continues on the face defined by $x_4=0$, alternating between green and red edges.  On the face $x_1=0$, the chains alternate between red and blue edges, except for the longest blue segment along $00cd.$ See figure 8 for $L'(3, 3)$. The simplex $L'(1, 3)$ is the chain $2001 \gtrdot 1101 \gtrdot 1010 \gtrdot 1002.$    

For $L'(5,3)$, we repeat the algorithm.  Figure 9 shows the faces for vanishing $x_1$ or $x_4.$  The remainder is filled out with $L'(3, 3)$ as noted.
\vspace{5pt}

The even case is handled in a similar manner, but with two layers for remaining chains.  We first add a red chain along the intersection of the faces on the outer layer. We continue adding chains as before, but shift to the inner layer to avoid the red chain.  Then $L(2t-4, 3)$ is inserted symmetrically by shifting multiplicities as before with the weak composition 2002.

See figures 10 through 13 for the case of $L'(4, 3)$.  In this case, the insertion of $L'(4, 0)$ is simply the addition of the node 2002.  For $L'(6, 3),$ we insert the simplex $L'(2, 3)$ corresponding to figure 3 with the symmetric chain decomposition for $L(3, 2).$

\begin{sidewaysfigure}
\vspace{300pt}
\includegraphics[width=.31\textwidth]{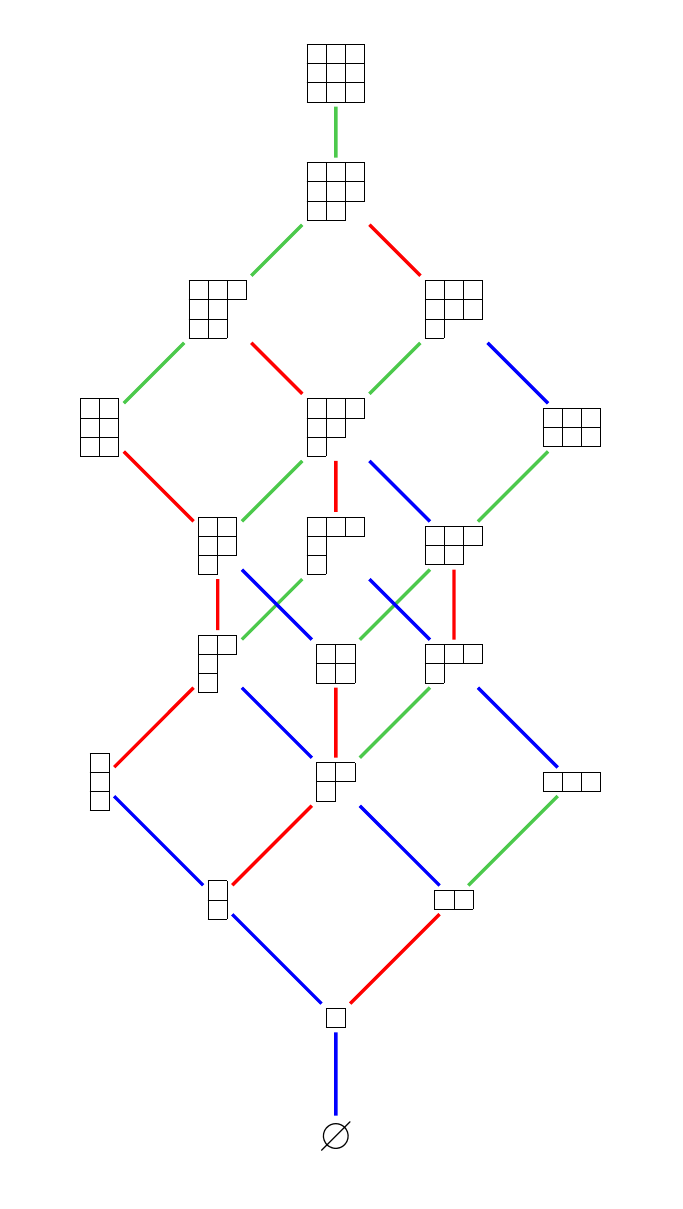}
\includegraphics[width=.35\textwidth]{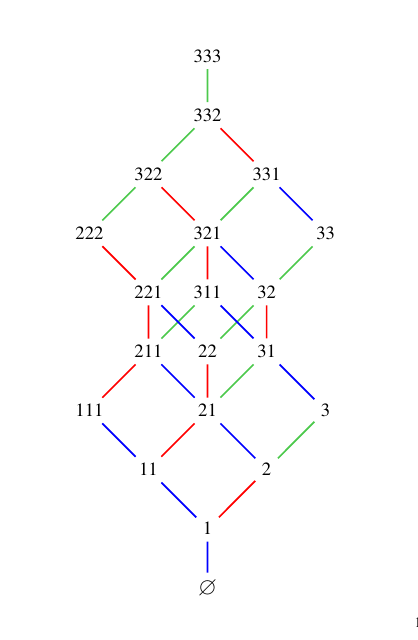}
\includegraphics[width=.31\textwidth]{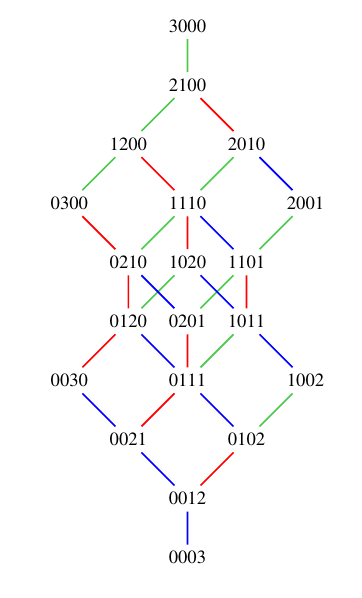}
\caption{Finite Young lattice $L(3,3)$ with colored links}
\label{}
\end{sidewaysfigure}

\begin{sidewaysfigure}
\vspace{300pt}
\hfil \includegraphics[width=.25\textwidth]{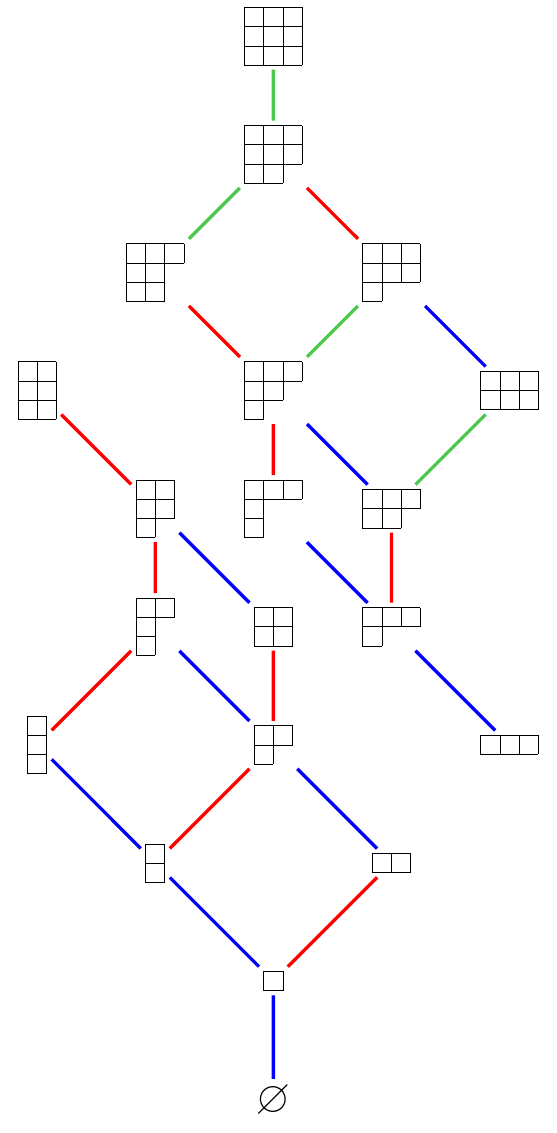} \hfil
\includegraphics[width=.26\textwidth]{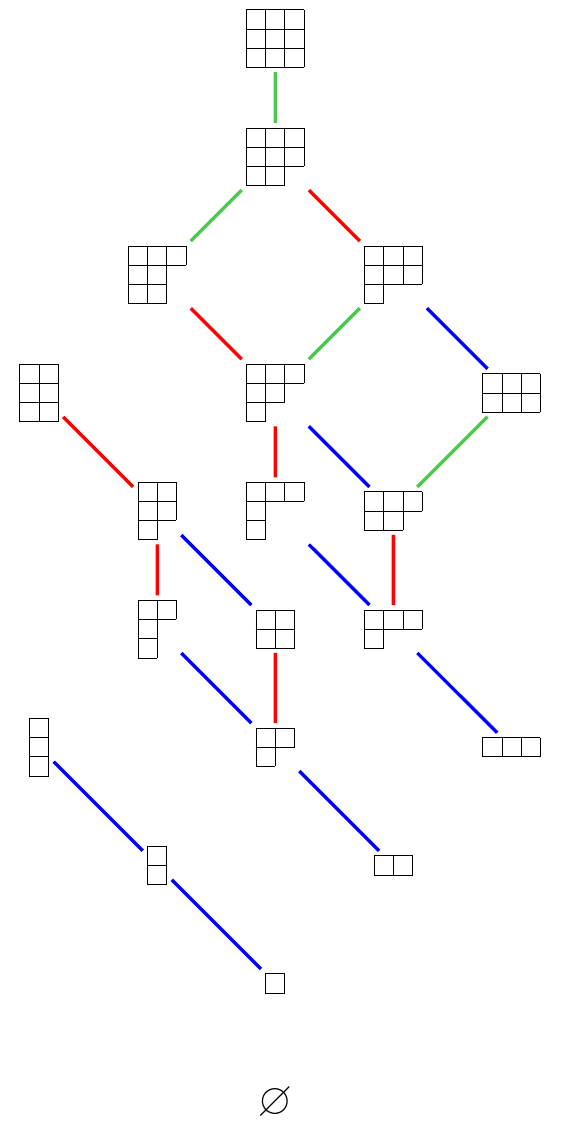}\hfil
\includegraphics[width=.25\textwidth]{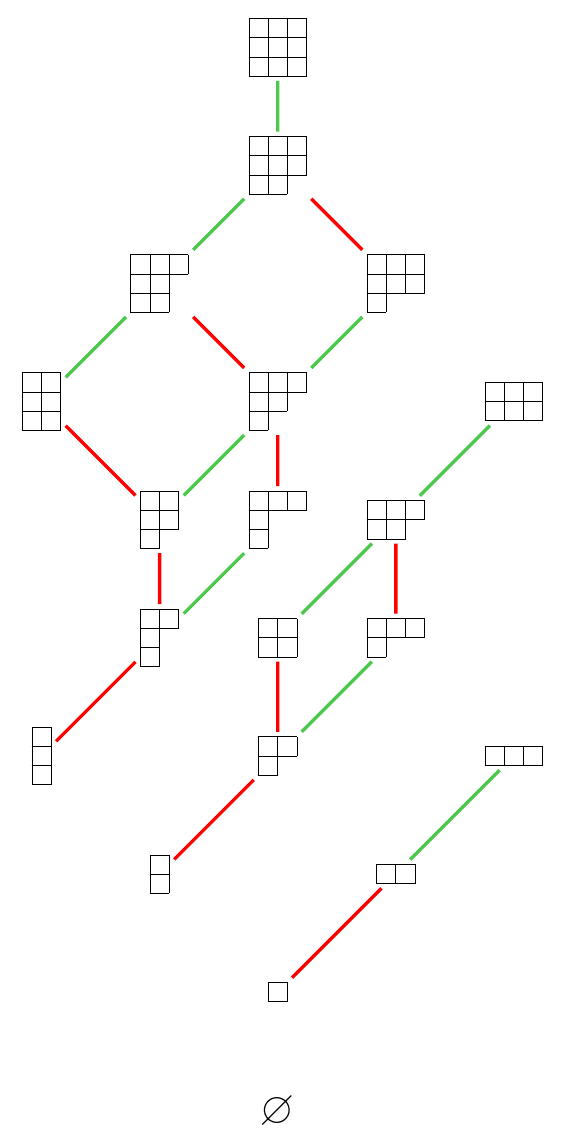}
\caption{Partial decompositions for $L(3,3)$}
\label{}
\end{sidewaysfigure}

\begin{figure}
\centering
\includegraphics[width=1.1\textwidth]{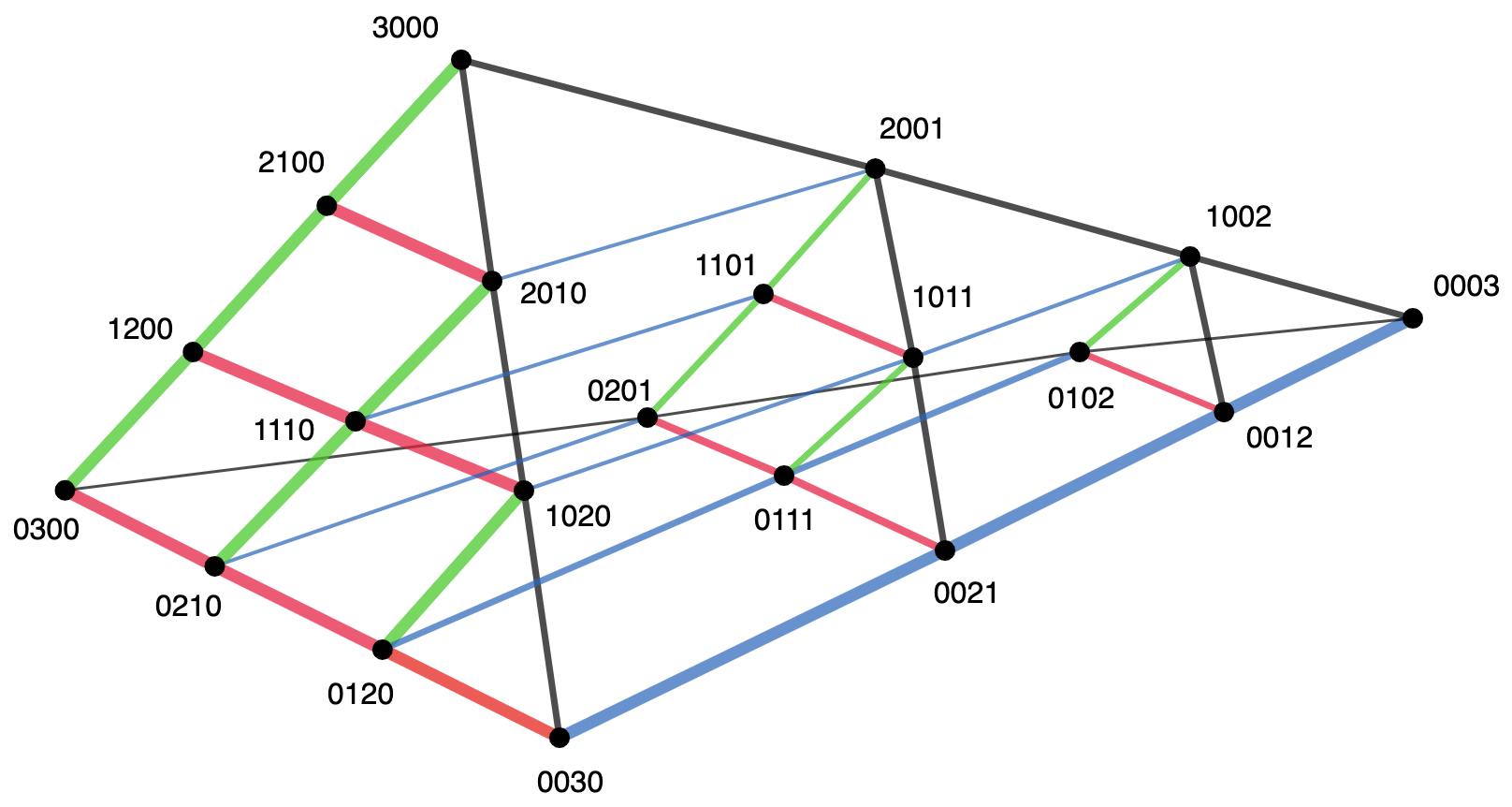}
\caption{ $L'(3, 3)$ with weight strings}
\label{}
\end{figure}
\vspace{20pt}

\begin{figure}
\centering
\includegraphics[width=1.1\textwidth]{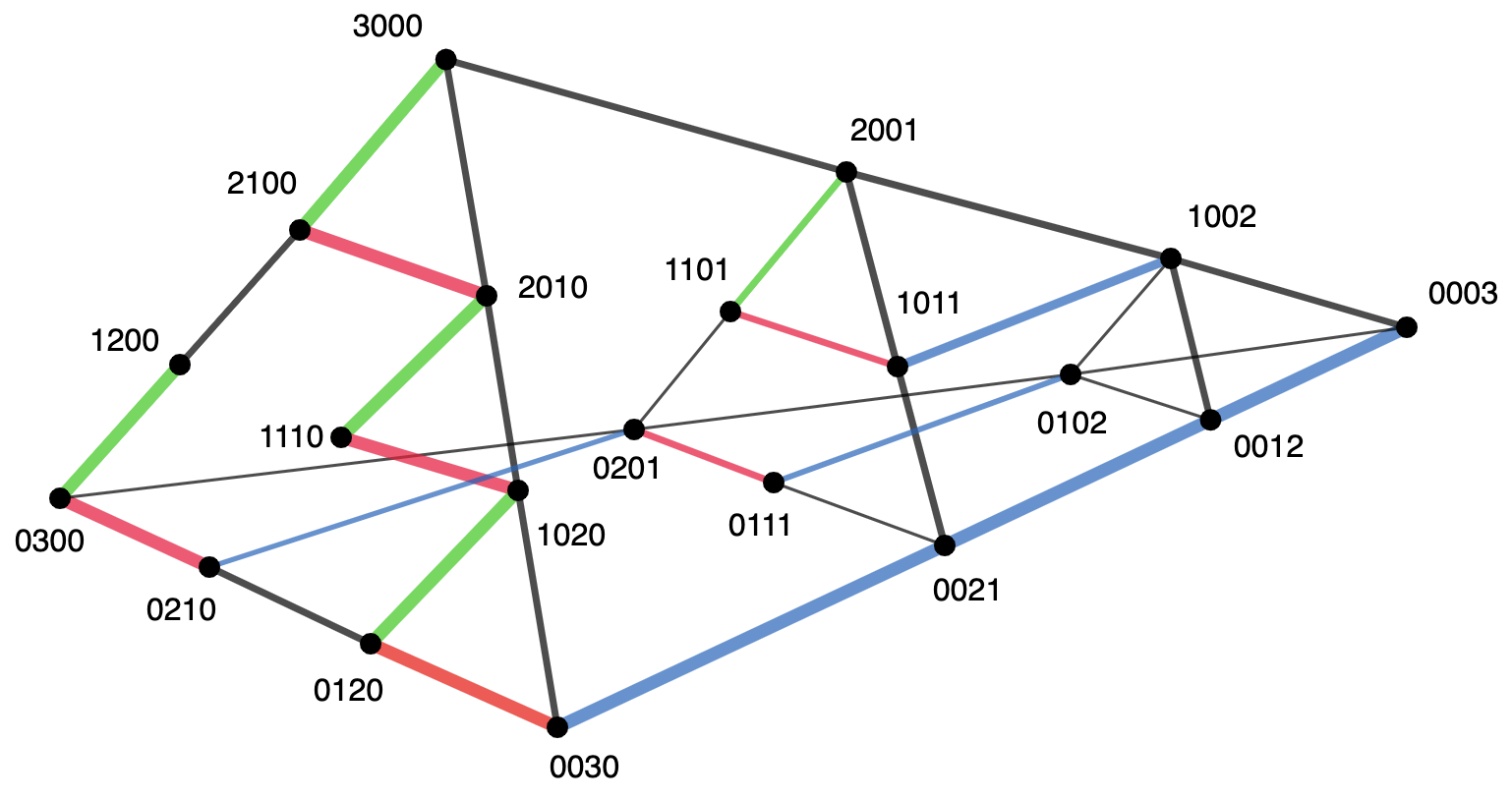}
\caption{Symmetric chain decomposition for $L'(3, 3)$ }
\label{}
\end{figure}

\begin{sidewaysfigure}
\vspace{350pt}
\includegraphics[width=1\textwidth]{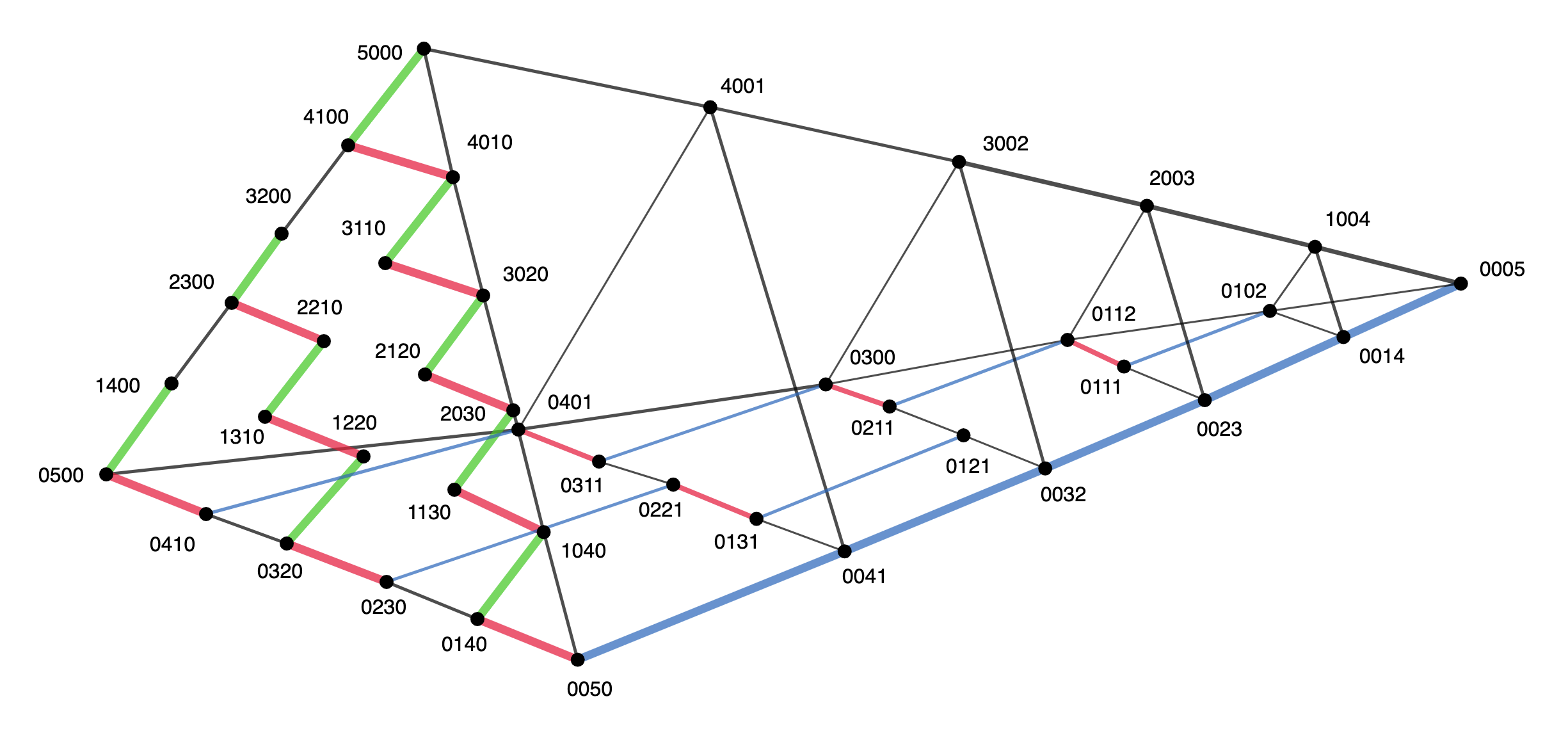}
\caption{Partial symmetric chain decomposition for $L'(5, 3)$ } \hfil
\end{sidewaysfigure}

\begin{figure}
\includegraphics[width=1\textwidth]{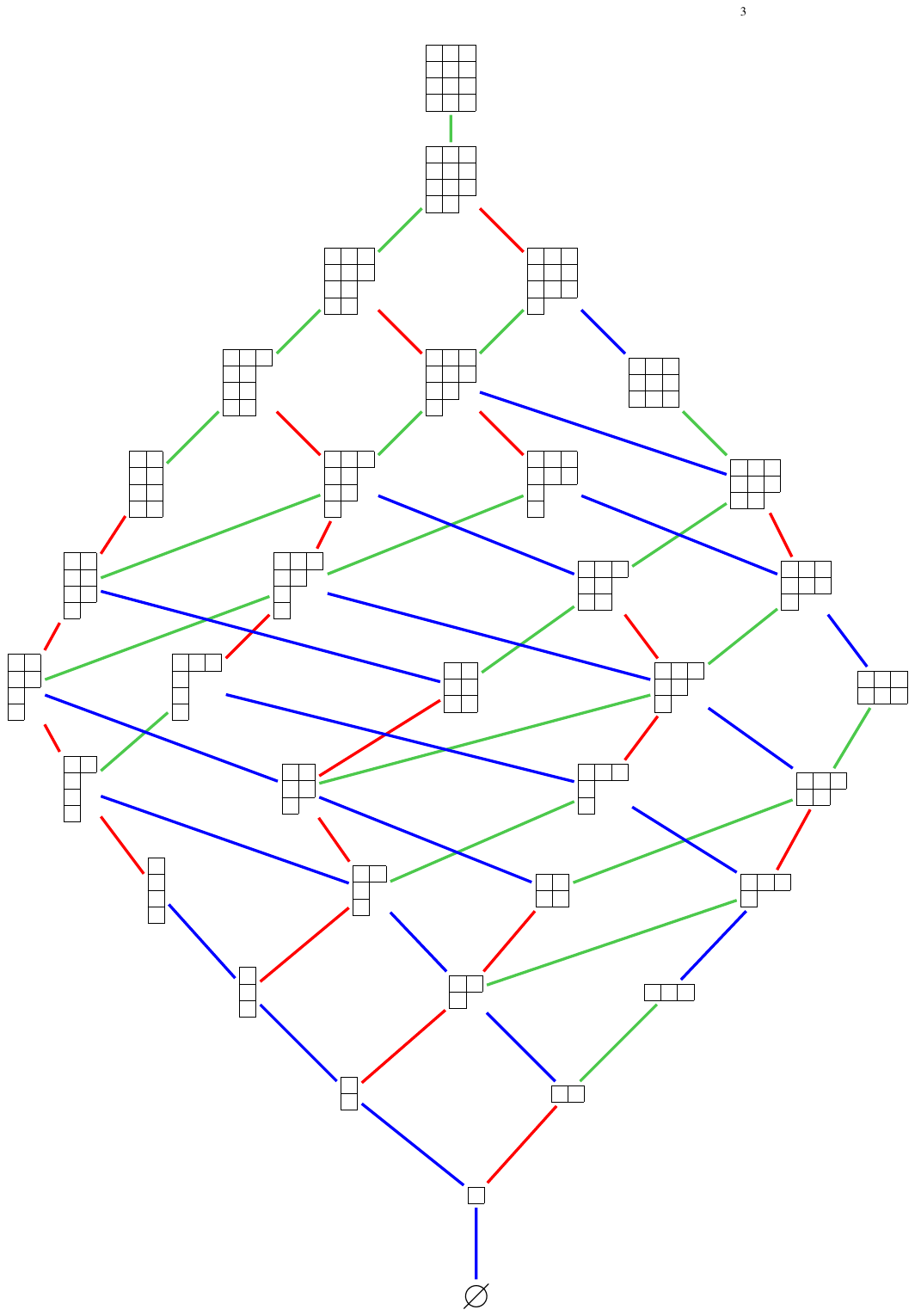}
\caption{FiniteYoung lattice $L(4,3)$ with colored links}
\label{}
\end{figure}

\begin{figure}
\includegraphics[width=1\textwidth]{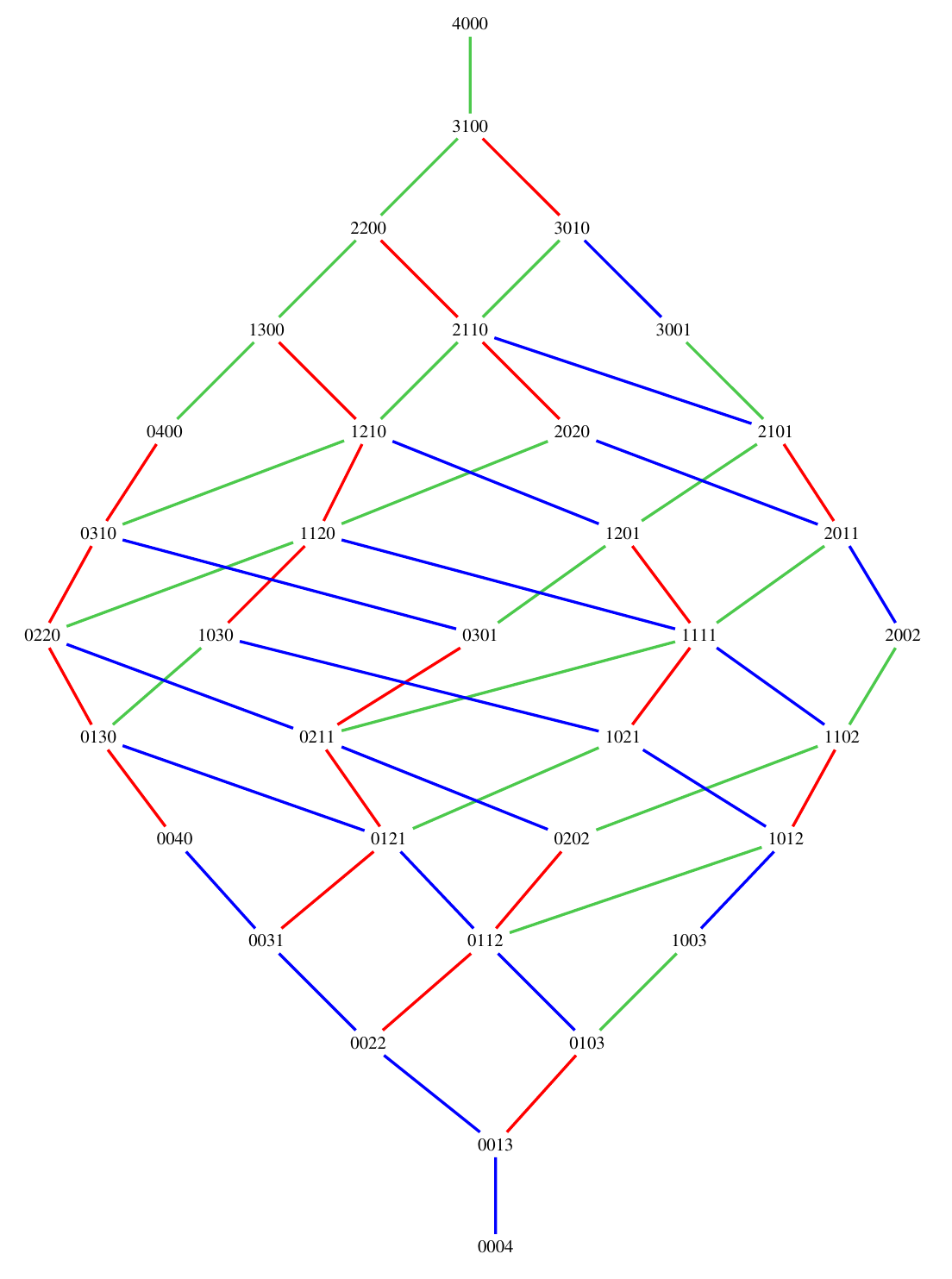}
\caption{ $L'(4,3)$ with multiplicities for nodes}
\label{}
\end{figure}

\begin{figure}
\includegraphics[width=1\textwidth]{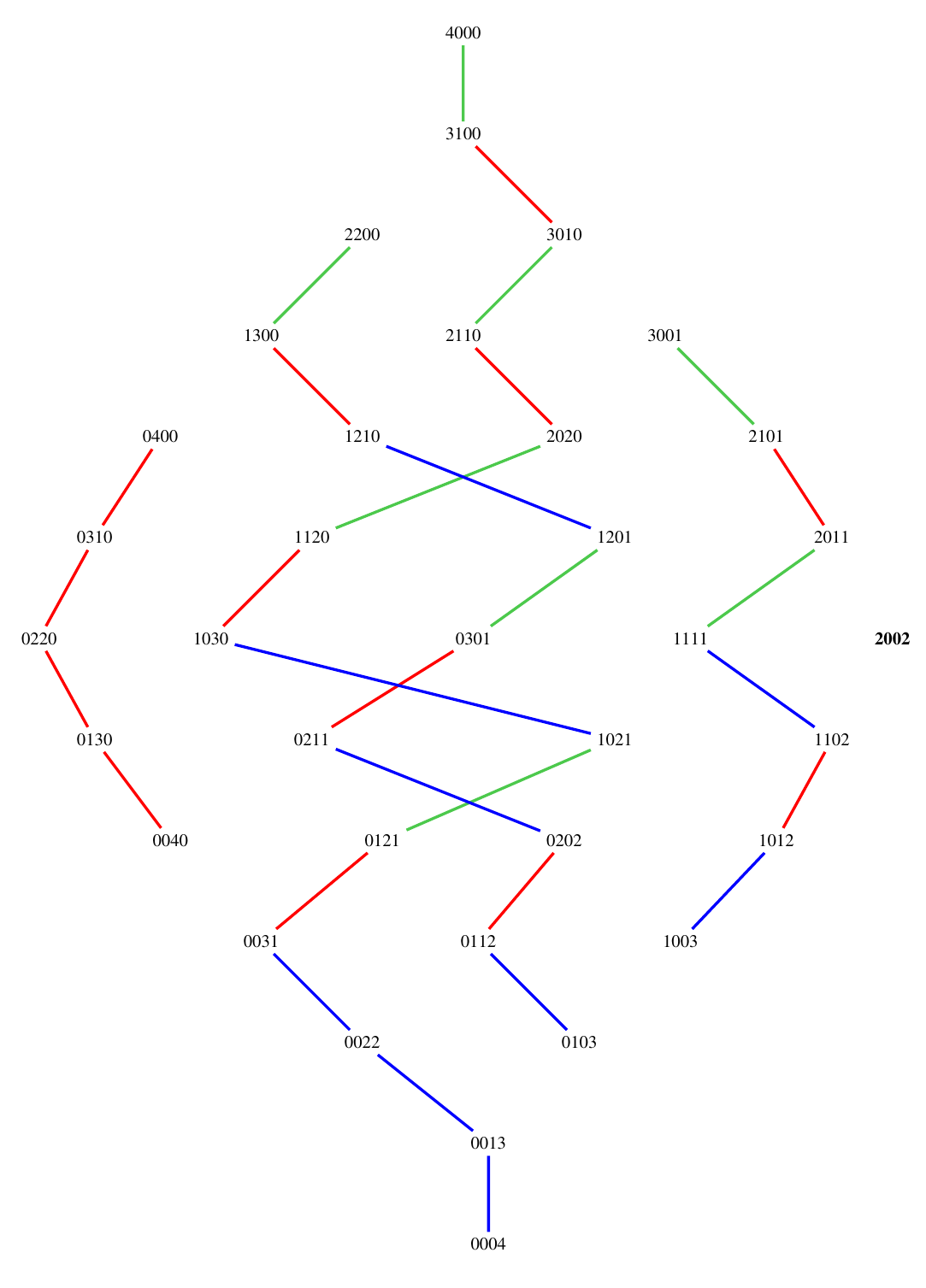}
\caption{Lindstr{\" o}m's algorithm for $L'(4, 3)$}
\label{}
\end{figure}

\begin{sidewaysfigure}
\vspace{275pt}
\includegraphics[width=1.1\textwidth]{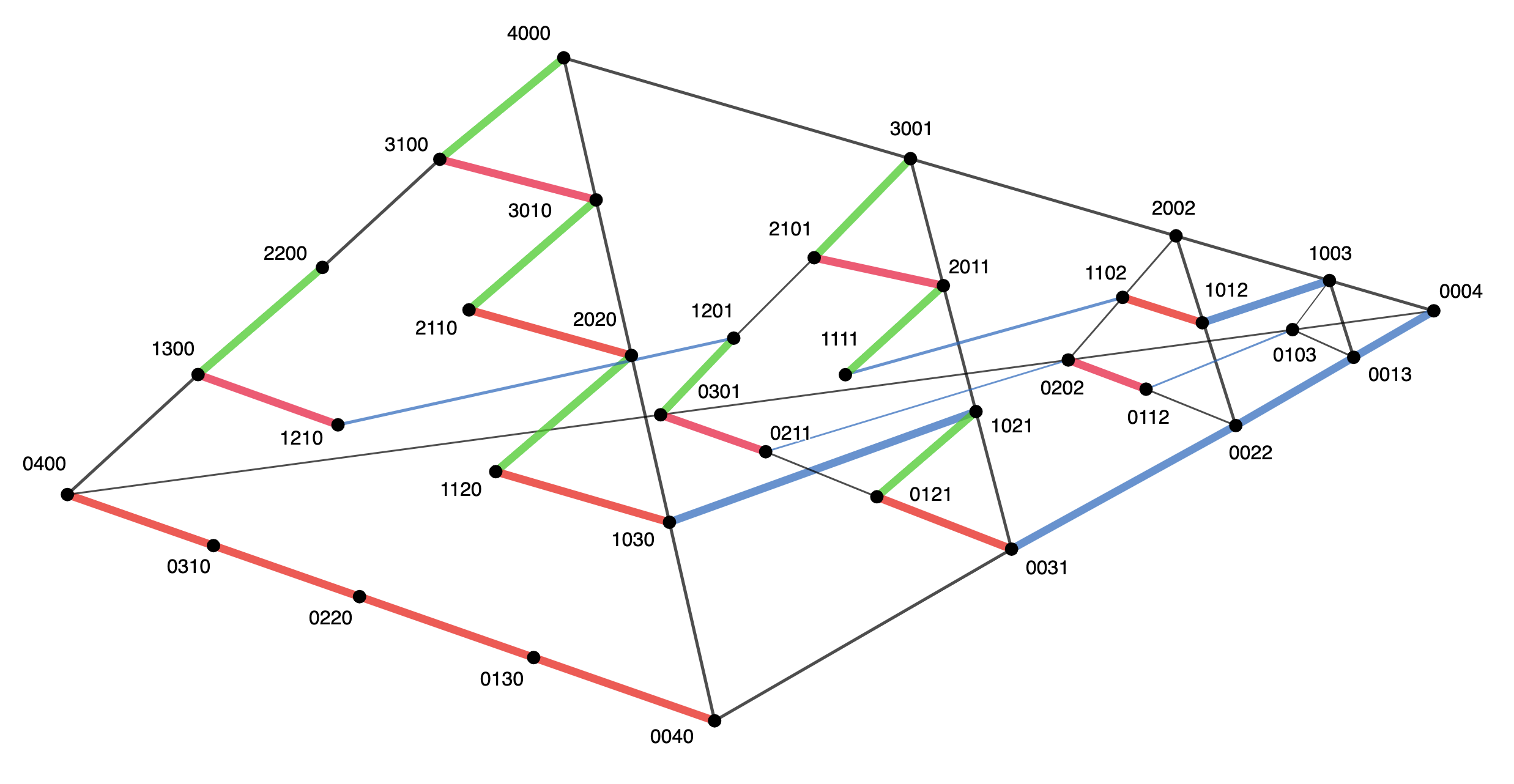}
\caption{Lindstr{\" o}m's algorithm for $L'(4, 3)$}
\label{}
\end{sidewaysfigure}

\bibliographystyle{amsplain}

\end{document}